\documentclass[12pt,a4paper,reqno]{amsart}
\usepackage{float}
\usepackage{amsfonts,amsmath,graphics,epsfig,latexsym,times}
\usepackage{ae,aecompl}
\usepackage[all]{xy}
\newif\ifpix \pixtrue

\usepackage[colorlinks=true]{hyperref}
\hypersetup{ pdfstartview=FitH }

\numberwithin{equation}{section}

\newcommand{\Mch}{M}
\newcommand{\Mtil}{\widetilde{M}}

\numberwithin{equation}{section}

\def\TT{{\mathbb{T}}} \def\PP{{\mathbb{P}}} \def\HH{{\mathbb{H}}}
\def\ZZ{{\mathbb{Z}}}  \def\QQ{{\mathbb{Q}}} \def\CC{{\mathbb{C}}}
\def\RR{{\mathbb{R}}} \def\LL{{\mathbb{L}}} \def\TT{{\mathbb{T}}} \def\EE{{\mathbb{E}}} 
\def\cO{{\mathcal{O}}}

\def\cK{{\mathcal{K}}}
\def\cV{{\mathcal{V}}}

\def\cU{{\mathcal{U}}}
\def\cR{{\mathcal{R}}}

\def\fM{{\mathfrak{M}}}
\def\fC{{\mathfrak{C}}}
\renewcommand{\epsilon}{\varepsilon}
\newcommand{\ip}[1]{\langle #1 \rangle}
\newcommand{\SU}{\mathrm{SU}}
\newcommand{\Hom}{\mathrm{Hom}}

\newcommand{\Ric}{\mathrm{Ric}}
\newcommand{\rank}{\mathrm{rank}  \;}

\newcommand{\CP}{\mathbb{CP}}

\newcommand{\del}{\partial}
\newcommand{\delb}{\bar\partial}

\newcommand{\ovS}{\overline{\Sigma}}
\newcommand{\chiorb}{\chi^{\mathrm{orb}}}

\newcommand{\geuc}{{g^{\mathrm{euc}}}}

\newcommand{\ovM}{\overline{M}}

\newcommand{\hM}{\widehat{M}}
\newcommand{\Mbar}{\overline{M}}
\newcommand{\Mhat}{\widehat{M}}
\newcommand{\Xbar}{\overline{X}}
\newcommand{\Xtil}{\widetilde{X}}
\newcommand{\Xhat}{\widehat{X}}
\newcommand{\Sbar}{\overline{\Sigma}}
\newcommand{\Scir}{\displaystyle {\mathop{\Sigma}^{\circ}}}
\newcommand{\Shat}{{\Sigma}}

\newcommand{\hS}{{\Sigma}}

\newcommand{\rd}{{\mathrm d}}

\newcommand{\oH}{\overline{\HH}}

\newtheorem{lemma}[subsubsection]{Lemma}
\newtheorem{prop}[subsubsection]{Proposition}
\newtheorem{cor}[subsubsection]{Corollary}
\newtheorem{theo}[subsubsection]{Theorem}
\newtheorem{theointro}{Theorem}

\newtheorem{propintro}[theointro]{Proposition}

\newtheorem{corintro}[theointro]{Corollary}

\theoremstyle{definition}
\newtheorem{dfn}[subsubsection]{Definition}
\newtheorem{rmk}[subsubsection]{Remark}

\newtheorem*{rmk*}{Remark}
\newtheorem*{rmks*}{Remarks}
\newtheorem{rmks}[subsection]{Remarks}

\newenvironment{remark*}{\begin{rmk*} --- \normalfont} { \end{rmk*} }
\newenvironment{remarks*}{\begin{rmks*} \begin{enumerate} \normalfont}
{\end{enumerate} \end{rmks*} } 
 
\title[CSCK surfaces and parabolic polystability]{Constant scalar
  curvature K\"ahler
  surfaces and parabolic polystability} 
\author{Yann Rollin}
\address{Yann Rollin, Imperial College, London, UK}
\email{rollin@imperial.ac.uk}
\author{Michael Singer}
\address{Michael Singer, University of
  Edinburgh, Sco\-t\-land} 
\email{m.singer@ed.ac.uk}
\thanks{First author supported by a University Research
  Fellowship of the Royal Society}

\begin{document}
{\Huge \sc \bf\maketitle}
\begin{abstract}
A complex ruled surface admits an iterated blow-up encoded by 
a parabolic structure with rational weights.
Under a condition of parabolic
stability, one can 
construct  a K\"ahler metric of constant scalar curvature on
the  blow-up according to \cite{RolSinX04}. We present a
generalization of this construction to the case 
of parabolically polystable ruled surfaces. Thus we can produce numerous
examples  of K\"ahler surfaces of constant scalar curvature with
circle or toric symmetry. 
\end{abstract}

\section{Introduction}

Let $M$ be a geometrically ruled surface of genus $g$. Thus $M$ is the
projectivization of some rank-2 holomorphic vector bundle $E\to
\Sigma$, where $\Sigma$ is a closed Riemann surface of genus $g$. Let
$\hM$ be an iterated blow-up of $M$, so there is a sequence of
holomorphic maps
$$
\hM =M_r \longrightarrow M_{r-1} 
\longrightarrow \cdots \longrightarrow M_1
\longrightarrow M_0=M
$$
and $M_j$ is the blow-up of $M_{j-1}$ at a point. The
existence of {\em scalar-flat K\"ahler} metrics\footnote{We shall
  use the acronym SFK for ``scalar-flat K\"ahler''} on such $\Mhat$ was
investigated especially by Claude LeBrun and his coworkers in the
1980's and 1990's using a combination of explicit constructions and
(singular) complex deformation theory.  New gluing techniques and
constructions were introduced by Arezzo and Pacard~\cite{ArePacX04,
ArePacX05} and the authors in \cite{RolSinX04,RolSin05} 
giving many examples of K\"ahler metrics of constant scalar
curvature\footnote{We shall use the acronym CSCK for ``constant scalar
  curvature K\"ahler''} on such surfaces.   It was noted by various
authors that the existence of CSCK metrics on $\hM$ appeared to be
closely related to the parabolic stability of the underlying bundle
$E\to\Sigma$ (the parabolic structure encodes the iterated blow-up in
a manner described carefully in \S\ref{secparab}).

In this paper we study the case corresponding to {\em strict}
parabolic polystability. In particular $E = L_1\oplus L_2$ is a
direct sum of line-bundles and $M$ contains two sections $S_j =
\PP(L_j)$. There is a holomorphic $\CC^*$ action preserving the fibres
and fixing $S_1$ and $S_2$, this action lifts to $\hM$, and our construction
gives many examples of CSCK and SFK metrics admitting an isometric
$S^1$-action.

We note in passing that a classification of SFK metrics with
non-trivial isometry group was claimed in Proposition~3.1 of
\cite{LebSin93}.  Unfortunately, as it was pointed out in
\cite{KimLebPon97}, there is an error at the end of the proof given
there and the conclusion that the genus of $\hM$ must be $\geq 2$ is
false. This paper gives scores of `counterexamples': SFK metrics with
non-trivial isometry group on blown-up ruled surfaces of genus $0$ and
$1$.  

\subsection{Statement of results}
The main result of this paper is the following theorem. 
\begin{theointro}
\label{maintheo}
Let $\Mch\to \hS $ be  a  parabolically polystable and non-sporadic   ruled
surface over a smooth Riemann surface with
 rational weights. Then its corresponding iterated blow-up $\Mhat$
carries a  CSCK metric. If $M\to\Sigma$ 
is strictly parabolically polystable,  $\Mhat$ admits furthemore a
non-trivial holomorphic 
vector field.

In addition, if the orbifold Riemann
surface $\Sbar$ deduced from $\Sigma$ and the parabolic structure
satisfies   $\chiorb(\Sbar)<0$, we may assume that the 
metric is SFK. 
\end{theointro}
The definitions of parabolic structures,  iterated 
blow-ups, and sporadic structures are all given in
\S\ref{secparab}. 
The definitions of  $\Sbar$ and its
orbifold Euler characteristic $\chiorb$ are recalled in \S\ref{secgood}.

\begin{rmk}
If the
condition ``parabolically polystable'' is replaced with the stronger
assumption ``parabolically stable'',
this result was already proved by the authors~\cite{RolSin05,
RolSinX04}.  So it is sufficient to prove Theorem~\ref{maintheo} 
 for strictly parabolically polystable ruled surfaces.
\end{rmk}

\begin{rmk} LeBrun's explicit construction of SFK metrics on blown-up
  ruled surfaces \cite{Leb91} gives rise to a SFK metric on a simple
  blow-up $\hM^s$, say, of a parabolically polystable ruled surface
  $M$. The surface $\hM$ is an iterated blow-up of $\hM^s$ also admits
  a SFK metric according to Theorem~\ref{maintheo}.
 However
  these metrics are {\em not} close to one another in a sense that
  will be  made
  more precise at ~\S\ref{sec:compare}.  Thus the parabolic structure can be used in
  different ways to encode CSCK metrics and it is not clear which is
  the most natural.
\end{rmk}
 In the parabolically stable case, it
was proved that in addition to the conclusion of
Theorem~\ref{maintheo},   any further blow-up of $\Mhat$ carries a
CSCK 
metric as well \cite{RolSinX04}.
This result extends  as follows:
\begin{propintro}
\label{prop:further}
  Let $\pi:\Mch\to \hS $ be a 
strictly parabolically polystable and non-sporadic  ruled
  surface with rational weights and $\Mhat$ be its iterated blow-up. 
Given any  finite collection of points $\{y_1,\cdots,
y_m\}\subset M\setminus \pi^{-1}(\{P_j\})$, where $P_j$ are
 the parabolic points of
$\Shat$, we define $\Mtil\to\Mhat$ by making further blow-ups at each $y_j$.

If the parabolic structure is not trivial and  $\Sigma$ is not the
sphere with exactly two parabolic points, then 
$\Mtil$ carries a CSCK metric.
If moreover $\chiorb(\Sbar)<0$, we may assume that the 
metric is SFK. 
\end{propintro}

\begin{rmk}
If we do not exclude   the trivial  parabolic structure or the case where $\Sigma$ is
the two-punctured sphere, the conclusion of Proposition~\ref{prop:further}
still holds, but
only for special configurations of points $y_1,\cdots, y_m$
(cf. \S\ref{sec:special}). 
\end{rmk}

 The proof of 
Theorem~\ref{maintheo} and Proposition~\ref{prop:further} relies on an extension of Arezzo-Pacard
gluing  theory~\cite{ArePacX05} for the orbifold setting
(cf. \S\ref{sec:ap}).  However, this extension is not
straightforward if the structure is sporadic and some significant work
would be needed to develop the appropriate gluing theorem in this case
(cf. Remark~\ref{rmk:laborious}).  We conjecture that 
Theorem~\ref{maintheo} also holds for sporadic parabolically polystable ruled surfaces.

The technical
difficulty is related to the following fact of independent interest:
for a pair of coprime integers $0<p< q$, consider the action of the cyclic subgroup $\Gamma_{p,q}\subset U(2)$
on $\CC^2$
generated by
\begin{equation}
\label{eq:action}
(z_1,z_2)\mapsto (z_1\zeta, z_2\zeta^{p}).
\end{equation}
where $\zeta$ is a $q$-th root of unity.
There is well known minimal resolutions $Y_{p,q}\to
\CC^2/\Gamma_{p,q}$ of the orbifold singularity called the
Hir\-ze\-bruch-Jung resolution and the resolution carries an
asymptotically locally Euclidean metric (ALE) which is moreover SFK by
a result of Calderbank-Singer~\cite{CalSin04} (particular cases are
due to Kronheimer, Joyce and LeBrun). Then we have the following result.
\begin{theointro}
\label{theo:mass}
  The ALE SFK metric constructed by Calderbank-Singer on $Y_{p,q}$ has
  non-positive mass. It has zero mass if and only if $p=q-1$, that is
  whenever $Y_{p,q}\to\CC^2/\Gamma_{p,q}$ is a crepant resolution.
\end{theointro}
\begin{rmk}
  This result generalizes a computation of LeBrun for the
  case $p=1$, thus giving a larger class of counterexamples to the
  \emph{generalized positive action conjecture} (cf.~\cite{Leb88}).
\end{rmk}

The general construction of CSCK metric given by
Theorem~\ref{maintheo} leads to numerous examples. Indeed, parabolic 
ruled surfaces are generically stable.
On the other hand we are interested in the following question: what
are the simplest SFK rational surfaces? 
This question was answered  in~\cite{RolSin05}  by showing that there
exists a SFK metric on certain $10$-points blow-ups of $\CP^2$ --- the minimal
number of times one has  to blow-up so that there is no obstruction for
existence a SFK metric. 
 By construction, these examples have no
non-trivial holomorphic vector fields. In the next proposition, we are
trying to answer the same question for metrics with an $S^1$-symmetry,
that is in 
presence of a non-trivial holomorphic vector field.  The case of a $15$-points blow-up found in 
  \cite{KimLebPon97} is improved as follows:
\begin{corintro}
\label{prop:simplest}  
The complex plane $\CP^2$ blown-up at $11$-suitably chosen points
  carries a non-trivial holomorphic vector field and a SFK metric.
\end{corintro}
\begin{rmk}
  Notice that it is not interesting to ask for more
  symmetries. Indeed, following an argument of Yau \cite{Yau74}, if a
  SFK surface $(X,\omega)$ has two linearly independent holomorphic
  vector fields   $X_1$, $X_2$, then $X_1\wedge X_2$ is a non trivial section of
  $K_X^{-1}$. However $4\pi c_1(X)\cdot [\omega]=\int s d\mu= 0$ and
  it follows that 
 $K_X$ is trivial. We deduce that $(X,\omega)$ is covered by a flat
 torus.
\end{rmk}
Now the question for  $10$-points remains open: it might be the case
that  
 $10$-points
blow-ups of $\CP^2$ with a non-trivial holomorphic vector field cannot
carry SFK metric. Yet there  is no obvious
obstruction known at the moment.

\subsection{A toric example}
\label{sec:simpex1}
 We will
illustrate  Theorem~\ref{maintheo} by highlighting a particularly
easy to describe  toric example. As we shall see, this example
 turns out to be the blow-up of a
strictly parabolically polystable ruled surface
(cf. \S\ref{sec:simpex2}). Therefore the application cannot be 
deduced from~\cite{RolSinX04}  and 
Theorem~\ref{maintheo} is required.

Consider the complex orbifold 
$$\Mbar=(\CP^1\times\CP^1)/\ZZ_2$$
 where
$\ZZ_2$ acts by  inversion $[z_0:z_1]\to [z_1:z_0]$ on each
factor and let $\Omega$ be an orbifold K\"ahler class on $\Mbar$, that
is an element of the orbifold De Rham cohomology $\Omega\in H^2_{\mathrm{DR}}(\Mbar,\RR)$ which can be
represented by an orbifold K\"ahler class. 
 Let $\Mhat$ be the  resolution of $\Mbar$ obtained by
replacing each of the four orbifold singularities $\CC^2/\ZZ_2$ of
$\Mbar$ with its resolution $\cO(-2)$. 

We will call
$E_1,\cdots, E_4$ the four exceptional divisors of self-intersection
$-2$.  Then we have the following result.
\begin{corintro}
\label{propexample}
  Let $\pi:\Mhat\to\Mbar$ be the resolution of
  $\Mbar=(\CP^1\times\CP^1)/\ZZ_2$ defined above, $\Omega$ be an 
 orbifold  K\"ahler 
  class on $\Mbar$ and  $\|\cdot\|$ be any fixed norm on $H^2(\Mhat,\RR)$. 
Then for all 
  $\epsilon >0$  there is a CSCK metric $\hat\omega$ on $\Mhat$ with 
\begin{equation}
\label{eq:kahlclass}
\|[\hat \omega] - \pi^* \Omega \|\leq \epsilon.
\end{equation}
\end{corintro}

\begin{rmk}
Notice that $\Mhat$ has two linearly independent holomorphic vector
fields. Therefore the CSCK metrics on $\Mhat$ have necessarily toric
symmetry. Many other CSCK metrics with torus symmetric can be constructed thanks to
Theorem~\ref{maintheo} (cf. Remark~\ref{rmk:toric}).
\end{rmk}

\begin{rmk}
No information is given about the K\"ahler class in Theorem
  \ref{maintheo}, however the estimate~\eqref{eq:kahlclass} is an
  immediate consequence 
of  the gluing theory used to prove the theorem (cf. Theorem~\ref{theoglue}).
\end{rmk}


\subsection*{Acknowledgments}
We thank Frank Pacard for making useful comments on
Theorem~\ref{theoglue}. We also wish to thank
Claude LeBrun for encouragements and for pointing out the remark made
in \cite{KimLebPon97} correcting \cite[Proposition 3.1]{LebSin93} (cf. \S\ref{sec:lsvs}).

\section{Parabolic structure on minimal ruled surfaces}
\label{secparab}
 
\subsection{Definitions}
A {\em geometrically ruled surface} $\Mch$ is by definition a minimal complex
surface obtained as $\Mch= \PP(E)$, where $E\rightarrow \Sigma
$ is a holomorphic vector bundle of rank $2$ over a \emph{smooth} Riemann surface 
$\Sigma$.
and we have an induced map $\pi:\Mch \rightarrow  \Sigma$ called the
{\em ruling}.

A parabolic structure on $\Mch$ consists of the following data:
\begin{itemize}
\item A finite set of distinct points $P_1,P_2,\cdots,P_n$ in
  $\Sigma$;
\item for each $j$, a choice of point  $Q_j \in F_j = \pi^{-1}(P_j)$;

\item for each $j$, a choice of  weight $\alpha_j \in   (0,1)\cap\QQ$.
\end{itemize}
A geometrically ruled surface with a parabolic structure will be called a
{\em parabolic ruled surface}.

If $S\subset\Mch$ is a holomorphic section of $\pi$, we define its
\emph{slope} by
$$\mu(S) = [S]^2 +\sum_{Q_j\not\in S}\alpha_j - \sum_{Q_j\in S}\alpha_j,
$$
where $[S]\in H_2(M,\ZZ)$ is the homology class of $S$ and
$[S]^2=[S]\cdot[S]$ is its self-intersection.

We say that a parabolic ruled surface is \emph{stable} (\emph{resp.}
\emph{semi-stable}) we have $\mu(S) >0$ (\emph{resp.}
$\mu(S)\geq 0$) for all holomorphic sections. 
We say that a parabolic ruled surface is \emph{polystable}  if it is
either stable, or semi-stable with two non-intersecting holomorphic 
sections $S_1$ and $S_2$ of  slope zero. We say that parabolic ruled surface is
\emph{strictly polystable} if it is polystable but not stable.

Alternatively, the expression ``$\Mch\to\Sigma$ is a
parabolically 
(poly)stable ruled surface'' means that the ruled surface is  endowed
with an \emph{a priori} fixed 
parabolic structure, and that it is (poly)stable w.r.t. that parabolic structure.
\begin{rmk}
For strictly parabolically polystable ruled surfaces,
 all  marked points $Q_j$ must lie either on $S_1$ or $S_2$.
\end{rmk}

 For some technical reasons, we will have to
exclude the so called \emph{sporadic} parabolic structures according to the
following definition.
\begin{dfn}[Sporadic parabolic structures]
\label{dfn:sporadic}
  Let $\Mch\to\Sigma$ be a parabolic ruled surface.
We will say that it is sporadic
if 
\begin{itemize}
\item it is strictly polystable
\item  $\Sigma$ is not the
sphere with exactly two marked points
\item every  parabolic point   $Q_j\in S_1$ has
 weight of the form $\alpha_j=\frac 1{q_j}$  and every
 parabolic point   $Q_j\in S_2$   has weight of the form $\alpha_j=\frac
 {q_j-1}{q_j}$ for some integer $q_j\geq 2$,
\item or, the same as above occurs with $S_1$ and $S_2$ exchanged.  
\end{itemize}
\end{dfn}

\subsubsection{Our toric example}
\label{sec:simpex2}
Consider the ruled surface $\Mch = \CP^1\times\CP^1$ with map
$\Mch\to\CP^1$ given by the first projection. Then, we introduce the
parabolic structure $P_1=[1:0]$, $Q_1=([1:0],[1:0])$, $P_2=[0:1]$,
$Q_2=([0:1],[0:1])$ with weights $\alpha_1=\alpha_2=\frac 12$.
\begin{lemma}
\label{lemma:ex}
  The ruled surface $\CP^1\times\CP^1\to\CP^1$ with the above parabolic
  structure is strictly parabolically polystable.
\end{lemma}
\begin{proof}
  Consider a constant section $S_1$ meeting $Q_1$ and $S_2$ meeting
  $Q_2$. They obviously do no intersect and $\mu(S_1)=\mu(S_2)=0$ (in
  particular the ruled surface is not parabolically stable).
If $S$ is any other constant section which does not meet any marked
point, then $\mu(S)=1$. If $S$ is a non-constant holomorphic section, $[S]^2\geq 2$ hence
$\mu(S)>0$. It follows from the discussion that the ruled surface is
not parabolically stable, but only polystable. 
\end{proof}
\begin{rmk}
  \label{rmk:toric0}
We obtain another example of stricly parabolically polystable ruled
surface 
by replacing the weights $\frac 12$ in the above example with any
rational weight such that
$\alpha_1=\alpha_2\in(0,1)$.
\end{rmk}
 
\subsection{Iterated blow-up of a parabolic ruled surface}
\label{secitbup}

Let $\Mch$ be a parabolic ruled surface. We shall now define a
multiple blow-up $\Mhat \to \Mch$ which is canonically
determined by the parabolic structure of $\Mch$.

In order to simplify the notation, suppose that the parabolic
structure on $\Mch$ is reduced to a single point $P\in \Sigma$; let $Q$ be the corresponding point in $F = \pi^{-1}(P)$ and let
$\alpha = \frac pq$ be the weight, where $p$ and $q$ are two coprime
integers such that $0<p<q$. Denote the Hirzebruch-Jung continued fraction
expansion of $\alpha$ by
\begin{equation}\label{e1.844}
\frac pq =  \cfrac{1}{e_1-\cfrac{1}{e_2-\cdots\cfrac{1}{e_l}}};
\end{equation}
define also
\begin{equation}\label{e20.844}
\frac {q-p}{q} = \cfrac{1}{e'_1-\cfrac{1}{e'_2-\cdots
\cfrac{1}{e'_m}}}.
\end{equation}
These expansions are unique if, as we shall assume, the $e_j$ and $e'_j$ 
are all $\geq 2$.

We give here a construction of the iterated blow-up $\hM$: the fiber $F$
has self-intersection $0$. The first step is to
blow up $Q$, to get a diagram of the form
\begin{equation}\label{e2.844}
\xymatrix{
{}\ar@{-}[rr]^{-1} && *+[o][F-]{}
\ar@{-}[rr]^{-1} &&{}
}
\end{equation}
By blowing up the intersection point of these two curves we get the
diagram
\begin{equation}\label{e3.844}
\xymatrix{
{}\ar@{-}[rr]^{-2} && *+[o][F-]{} 
\ar@{-}[rr]^{-1} && *+[o][F]{}
\ar@{-}[rr]^{-2} &&{} 
}
\end{equation}
Then we perform an iterated  blow-ups of one of the two intersection of the
only $-1$-curve in the diagram. Given $\alpha=p/q$, there is a unique way
(cf.~\cite[Proposition 2.1.1]{RolSin05}) to choose at
each step which
point has to be blown-up  in order to get the following diagram
\begin{equation}
\label{eq:mbup}
\xymatrix{
\ar@{-}[r]^{-e_1} & *+[o][F-]{}
\ar@{-}[r]^{ -e_{2}} &  *+[o][F-]{}
\ar@{--}[r] &  *+[o][F-]{}
\ar@{-}[r]^{-e_{l-1}} &  *+[o][F-]{}
\ar@{-}[r]^{-e_l} &  *+[o][F-]{}
\ar@{-}[r]^{-1} &  *+[o][F-]{}
\ar@{-}[r]^{-e'_{m}} &  *+[o][F-]{}
\ar@{-}[r]^{-e'_{m-1}} &  *+[o][F-]{}
\ar@{--}[r] &  *+[o][F-]{}
\ar@{-}[r]^{-e'_{2}} &  *+[o][F-]{}
\ar@{-}[r]^{ -e'_{1}} & 
}
\end{equation}
where the $-e_1$-curve is the proper transform of the original fiber
$F$.

More generally, if $\Mch$ has more parabolic points, we perform the same operation
for every point and get a corresponding iterated blow-up
$\hM\to\Mch$.

\subsubsection{Back to the toric example}
\label{sec:back}
Consider the blow-up $\Mhat\to\Mch$ of the parabolic ruled surface defined in
\S\ref{sec:simpex2}. The fiber over each marked point $P_j$, gives a
configuration of curves shown in \eqref{e3.844}. Now contract the
four $-2$ curves. Thus, we obtain a complex orbifold surface $\Mbar$
which is precisely $(\CP^1\times\CP^1)/\ZZ_2$ as described in
\S\ref{sec:simpex1} and $\Mhat$ is its resolution. 

\section{Representations and K\"ahler orbifolds}
\label{secfactory}
We review the construction of  K\"ahler orbifold of constant scalar
curvature
of~\cite{RolSin05}, \cite{RolSinX04}.

\subsection{Orbifold Riemann surfaces}
\label{secgood}
We start with a closed Riemann surface $\overline 
\Sigma $ of genus $g$ with a 
finite set of orbifold points $P_1,P_2,\cdots, P_k$, with local
ramified cover of order
$q_1,q_2,\cdots, q_k >1$.
 Recall first the description of the fundamental group of the 
punctured Riemann surface $ \Scir=\Sbar\setminus\{P_j\}$:
\begin{multline*}\label{e3.10.12.03}
\pi_1(\Scir)= \\
\langle a_1,b_1,\ldots,a_g,b_g,l_1,\ldots, l_k\; | \;
[a_1,b_1][a_2,b_2]\ldots[a_g,b_g]l_1\ldots l_k = 1\rangle
\end{multline*}
Here the $a_j$ and $b_j$ are standard generators of $\pi_1(\Sbar)$
 and  $l_j$ is
(the homotopy class of) a small loop around $P_j$. The orbifold
fundamental group is defined by
$\pi^\mathrm{orb}_1(\ovS) = \pi_1(\Scir) / G$,
where $G$ is the normal subgroup  of $\pi_1(\Scir)$ generated by  $l_1^{q_1}, l_2^{q_2},\cdots,
 l_k^{q_k}$.

The orbifold Euler characteristic is defined by
$$ 
\chi^\mathrm{orb}(\Sbar):= \chi^\mathrm{top}(\ovS)  - \sum_{j=1}^k
(1- \frac 1 {q_j}) .  $$ 

 Let us call an orbifold Riemann surface ``good''
if its orbifold universal cover admits a compatible K\"ahler metric of
CSC $\kappa_1$, say.  
However, not every orbifold Riemann
surface is ``good'' as explained in the next section.
\subsection{Facts about good Riemann surfaces}
The only orbifold Riemann surfaces which are not good are the one
topologically equivalent to $S^2$ with exactly one orbifold point (called
the  {\em tear-drop}) or with exactly two orbifold points of
distinct orders.

The following summarizes basic facts on good orbifold Riemann
surfaces: 
if $\ovS$ is good, then the sign of $\kappa_1$ is the same as
the sign of $\chi^{\mathrm{orb}}(\ovS)$ (by the Gauss-Bonnet theorem
for orbifolds). 

Recall that for  a complex manifold $X$, possibly with orbifold singularities, we
 will denote $\cV_0(X)$ 
the space of $(1,0)$-holomorphic vector fields vanishing
at some point on $X$.  

Orbifold ruled surfaces have no non-trivial holomorphic vector fields
except in the following cases
\begin{itemize}
\item $\Sbar =\CP^1$ or $\TT$ (a smooth elliptic curve),
\item $\Sbar = \CP^1/\ZZ_q$ for $q\geq 2$.
\end{itemize}
For more details, the reader can refer to \cite[\S2.1]{RolSinX04}.

\subsection{K\"ahler orbifold ruled surfaces}
\label{sec:korb}
In this section, we assume that $\ovS$ is a good orbifold Riemann
surface and we endow $\ovS$ with an orbifold K\"ahler metric $\bar
g^\mathrm{\Sigma}$ of constant curvature $\kappa_1$.  Note that, just
as for ordinary Riemann surfaces, we have $$
\ovS = \cU/\pi_1^{\mathrm{orb}}(\ovS)
$$ 
where the fundamental group acts by isometries on the universal
cover $\cU$ which is equal to $\HH$, $\EE$ or $\CP^1$, according as
$\kappa_1$ is negative, zero, or positive.

Let $g^\mathrm{FS}$ be the Fubini-Study metric with curvature
$\kappa_2>0$, on $\CP^1$.  The product metric is a  K\"ahler metric
of constant scalar curvature $s=2(\kappa_1+\kappa_2)$. Notice that
whenever $\chiorb(\Sbar)<0$, we have $\cU =\HH$. Hence we can choose
$\kappa_1=-\kappa_2$ and the metric is scalar-flat.

 We define the space of representations
$$
\cR(\Sbar) = \left \{ \rho \in \Hom(\pi_1^\mathrm{orb}(\overline\Sigma),
\SU(2)/\ZZ_2) \; |\; \rho(l_j) \mbox{ has order } q_j \right \}.
$$
For a given $\rho\in\cR(\Sbar)$, we deduce a faithful and isometric
action of $\pi_1^\mathrm{orb}(\ovS)$  on $\cU\times\CP^1$. 
The orbifold quotient
$$\overline M_\rho = \ovS\times_\rho\CP^1 =
(\cU\times\CP^1)/ \pi_1^\mathrm{orb}(\overline\Sigma)$$
inherits a CSCK metric denoted $\bar g_\rho$. Moreover, we may
assume that this metric is scalar-flat when $\chiorb(\Sbar)<0$.

As we shall see, the space of $(1,0)$-holomorphic vector fields
$\cV_0(\Mbar_\rho)$ on $\Mbar_\rho$ plays an essential role in the
gluing theory. The next proposition gives its dimension depending on $\rho$.
\begin{prop}
\label{prop:vf}
  Given  a good orbifold with no non-trivial holomorphic vector
  field $\Sbar$ and $\rho\in\cR(\Sbar)$, we have either:
  \begin{enumerate}
  \item \label{enum:0} no point in $\CP^1$ is fixed under the action of
    $\pi_1^\mathrm{orb}(\Sbar)$ via  $\rho$,
\item \label{enum:1} exactly two points in $\CP^1$ are fixed by group
  action,
\item \label{enum:2} or, $\rho$ is trivial.
  \end{enumerate}
In each case, we have
\begin{enumerate}
\item $\cV_0(\Mbar_\rho)$ is trivial,
\item $\dim_\CC \cV_0(\Mbar_\rho)= 1$,
\item or, $\dim_\CC \cV_0(\Mbar_\rho)= 3$.
\end{enumerate}
\end{prop}
\begin{proof}
  First notice that we are clearly in  one of the cases
  \ref{enum:0}--\ref{enum:2}. Indeed, either $\rho(\pi_1^{\mathrm{orb}}(\Sbar))$ contains two rotations with
  distinct axes, and we are in case \ref{enum:0},  or it consists only of rotations about a common axis
  and we are either in case \ref{enum:1} or \ref{enum:2}. 

The fact that $\cV_0(\Mbar_\rho)$ is trivial in case \ref{enum:0} is
proved in \cite[Theorem 3.4.1]{RolSinX04}.  By reading the proof, we
see that, more generally, real holomorphic vector fields, correspond to
an infinitesimal isometry, which is given in this case,  up to scaling
by a constant, by a
pair of antipodal points on $\CP^1$, fixed by $\rho$. The conclusion on
the dimension follows.
\end{proof}

We will also need the following result:
\begin{prop}
\label{prop:vf2}
Assume that 
  $\Sbar =\CP^1/\ZZ_q$ for some integer $q>1$ and let
  $\rho\in\cR(\Sbar)$.
Then $\dim_\CC \cV_0(\Mbar_\rho)= 2$. However $\Mbar_\rho$ can be
endowed with an isometric $\ZZ_2$-action such that there are no
non-trivial $\ZZ_2$-invariant holomorphic vector fields. 
\end{prop}
\begin{proof}
 We consider $\CP^1$ with its standard Fubini-Study metric and the
 (isometric) inversion
$$
[z_0:z_1]\mapsto[z_1:z_0].
$$

Suppose that we are given a rotation of  $\CP^1$ the form
$$
[z_0:z_1]\mapsto [\zeta z_0:z_1],
$$
where $\zeta$ is a $q$-th root of unity.
 The space of
holomorphic vector fields on
$\Sbar$ is $1$-dimensional and spanned by $\delb^\sharp\phi$
 where
\begin{equation}
  \label{eq:pot}
\phi= \frac {|z_0|^2-|z_1|^2}{|z_0|^2+|z_1|^2}
\end{equation}
is a function with average $0$. It is readily checked
that $\delb^\sharp\phi$ is not $\ZZ_2$-invariant.

The isometric inversion of $\CP^1$ acts  diagonally on
$\CP^1\times\CP^1$
and descends to an isometric $\ZZ_2$-action on the quotient
$\Mbar_\rho=(\CP^1\times\CP^1)/\pi^\mathrm{orb}_1(\Mbar_\rho)$. 

It is not hard to check that the space of holomorphic vector fields on
$\Mbar_\rho$ is spanned by the vector fields with potential  $\phi_1$
and $\phi_2$ coming from each factor. However these vector fields are
not $\ZZ_2$-invariant.
\end{proof}

\begin{rmk}
  The case \ref{enum:0} of Proposition \ref{prop:vf} is precisely the one studied in
  \cite{RolSin05} and \cite{RolSinX04}. The case \ref{enum:2} leads to a
  trivial question. Indeed, having a trivial $\rho$
  implies that $\Sbar\simeq \Sigma$ has no orbifold points, and 
  $\Mbar_\rho\simeq \Sigma\times \CP^1$. Thus, $\Mbar_\rho$ is smooth
  and carries obvious CSCK metrics given by the product of
  metrics of constant curvature on each factor.
Therefore,  we shall focus on case \ref{enum:1} in this paper.
\end{rmk}

\subsection{Desingularization of orbifold ruled surfaces}
\label{sec:desors}
More concretely, near an orbifold point $P$ of order $q$, the Riemann surface
$\Sbar$ , is uniformized by $\Delta/\ZZ_q$, where $\Delta$ is a small
disc centered at the origin in $\CC$. Then $\pi^{-1}(\Delta/\ZZ_q)\subset\Mbar_\rho$ is
given by the quotient $(\Delta\times \CP^1)/\ZZ_{q}$ and the action of
$\ZZ_q$ is
generated by
$$
(z,[u:v])\mapsto (z\zeta, [u\zeta^{\frac p2}:v\zeta^{-\frac p2}])
$$
where $\zeta$ is a $q$-root of unity, $p$ is an integer coprime with
$q$ and such that $1\leq p <q$. 

Notice that there are two orbifold points $A$ and $B$ in $(\Delta\times
\CP^1)/\ZZ_{q}$ at the points $A=(0,[0:1])$ and $B=(0,[1:0])$. Using the affine
coordinate $v=1$, we see that the singularity near the  orbifold
point $A$ is modelled on
$\CC^2/\Gamma_{p,q}$, where the action $\Gamma_{p,q}$ is
generated by
$$
(z,u)\mapsto (z\zeta, u\zeta^{p}).
$$
The other orbifold singularity at $B$ is given similarly by $\CC^2/\Gamma_{q-p,q}$.
There are well known minimal resolutions $Y_{p,q}\to
\CC^2/\Gamma_{p,q}$ called Hir\-ze\-bruch-Jung resolutions. By gluing them
at each orbifold point, we get a resolution denoted
$$\Mhat_\rho\to\Mbar_\rho.$$

\subsection{The theorem of Mehta-Seshadri}
We unravel how the stability condition for a parabolic ruled
surface is related to the construction of a CSCK metric on its
blow-up. Given a
parabolic ruled surface $\Mch \to \hS$,
 we deduce a orbifold Riemann surface $\ovS$ by introducing an orbifold
 singularity of order $q_j$ at every parabolic point $P_j\in \hS$ of weight
 $\alpha_j=p_j/q_j$. As a corollary of Mehta-Seshadri theorem~\cite{MehSes80}, we have the
 following proposition.
\begin{prop}
\label{prop:MS}
Let $\Mch\to\hS$ be a parabolically polystable ruled surface with
rational weights. Then
there exists a homomorphism 
$\rho \in \cR(\Sbar)$
such that  $\Mhat \simeq
\Mhat_\rho$, where $\hM$ is the iterated blow-up of the parabolic ruled
surface $\Mch$  described in \S\ref{secitbup} and $\hM_\rho$ is the
 resolution of the orbifold $\ovM_\rho$ defined in 
 \S\ref{secfactory}. Furthermore, $\rho(l_j)$ is conjugate to 
$$
\pm \left ( 
\begin{array}{cc}
  e^{i\pi\alpha_j} & 0 \\
0 &   e^{-i\pi\alpha_j}
\end{array}
\right ).
$$
In addition, $M\to\Sigma$ is parabolically
 stable if and only if $\rho$ does not fix any point in $\CP^1$.
\end{prop}
\begin{proof}
This Proposition is just a reformulation of the Mehta-Seshadri
theorem, usually stated with the language of  holomorphic vector
bundles~\cite{MehSes80}.
 For more details about our point of view, the reader may refer to \cite[Theorem
 3.3.1]{RolSin05}.
\end{proof}

\section{Desingularization of K\"ahler orbifolds}
We move to a more general setting where $\Xbar$ is a compact complex
surface with isolated orbifold singularities $x_1,\cdots, x_k$. By
definition, an orbifold point, $\Xbar$ is uniformized by a neighborhood
of $0$ in $\CC^2/\Gamma_{p,q}$. 

Similarly to the case of orbifold ruled surfaces (cf. \S
\ref{sec:desors}) we can consider the minimal resolution 
$\Xhat\to \Xbar$
 obtained by
gluing in Hirzebruch-Jung resolutions $Y_{p,q}\to\CC^2/\Gamma_{p,q}$ at
orbifold points. 
Given a finite collection of points $\{y_1,\cdots,y_r\}\subset
\Xbar\setminus\{x_1,\cdots,x_k\} $ we define the blow-up
$\Xtil\to\Xhat$
 of $\Xhat$ at each $y_j$.

\subsection{Gluing for CSCK metrics}
\label{sec:ap}
In this section, we assume that $\Xbar$ carries a CSCK metric, in the orbifold sense.

\subsubsection{Asymptotics of the Calderbank-Singer ALE spaces}
\label{sec:calsin}
By a result of Calderbank-Singer \cite{CalSin04}, the
Hirzebruch-Jung resolution $Y_{p,q}$ carries  ALE scalar-flat 
K\"ahler metrics. The asymptotics of the Calderbank-Singer  metrics are
carefully investigated in
\S\ref{sec:asymptotics}:
on  the chart 
$(\CC^2\setminus 0)/\Gamma_{p,q}$ the metric can be written in the form
\begin{equation}
  \label{eq:asympt}
   \omega =   dd^c (|z|^2+ m \log |z|^2 + u),
\end{equation}
 for some constant $m\in \RR$ and  $u=\cO(|z|^{-1})$. Moreover,
\begin{enumerate}
\item if $0<p<q-1$ we have $m<0$,
\item if $0<p=q-1$ then $m=0$.
\end{enumerate}
Notice that the Burns metric on $\CC^2$ blown-up at the origin can be
written in the form~\eqref{eq:asympt} as well. However, we have $m>0$
in this case.

\subsubsection{The deformation theory}
 One can patch them with the orbifold CSCK
metric on $\Xbar$ in order to get  approximate smooth CSC
K\"ahler metrics on the resolution $\Xhat$. There is gluing theory based 
on this picture and developed by Michael Singer and the author for
scalar-flat K\"ahler 
surfaces in \cite{RolSin05} and, later, for CSCK manifolds
by Arezzo-Pacard~\cite{ArePacX04}.

 The idea builds on the deformation theory  for CSC
 K\"ahler metric studied  by Simanca-LeBrun~\cite{LebSim94}: we perturb the
approximate CSCK metric on $\Xhat$ and show,
using the implicit function theorem, that there is a nearby
CSCK metric.

The fourth order linear operator
$$
\LL_{\Xbar}= - \frac 12 \Delta_{\Xbar}^2 - \Ric_{\Xbar}\cdot \nabla_{\Xbar}^2
$$
plays a central role in the gluing theory, since it is  to the
linearization of the map $\phi\mapsto s_\phi$, where $s_\phi$ is the
scalar curvature of the K\"ahler metric
$\omega_\phi=\omega_{\Mbar}+i\del\delb\phi$. We will denote 
$$
\cK(\Xbar)=\left \{\phi\in\ker\LL_{\Xbar}\; \mbox{ s.t. }\;
  \int_{\Xbar}\phi=0 \right\}
$$
where $\phi$ are real valued functions.
In the case of a CSCK metric we have the formula
\begin{equation}
\label{eq:formvf}
\LL_{\Xbar}=(\delb\delb^\sharp)^*\delb\delb^\sharp.
\end{equation}
It follows from \eqref{eq:formvf} that 
$$\cK(\Xbar)\otimes\CC \simeq \cV_0(\Xbar)
$$
$\cV_0({\Xbar})$ 
where  $\cV_0(\Xbar)$ is the space of $(1,0)$-holomorphic vector fields which
vanish at some point in $\Xbar$.  Furthermore, the correspondence  is
given by the map
$$
\phi\mapsto \delb^\sharp\phi,
$$
 extended by $\CC$-linearity.
The operator $\LL_{\Xbar}$ is elliptic, hence $\cK(\Xbar)$ is finite
dimensional and its kernel is spanned a by
$\phi_0,\phi_1,\cdots,\phi_r$, where $\phi_0=1$ and
$\phi_j\in\cK(\Xbar)$ for $j\geq 1$ form a linearly independent family.

\subsubsection{Generalization Arezzo-Pacard gluing theory}
Let $x_1, \cdots, x_k$ be the orbifold singularities in $\Xbar$
modelled on $\CC^2/\Gamma_{p_j,q_j}$  where $(p_j,q_j)$ are coprime integers with
$0<p_j<q_j$. We arrange so that for some $l\leq k$, the points
$\{x_1,\cdots,x_l\}\subset\{x_1,\cdots,x_k\}$ is exactly the subset of
points for which $p_j\neq q_j-1$. 
Let
$y_1,\cdots,y_m$ be a collection of smooth points and $\Xbar$.
 We introduce
the matrix
\begin{equation}
\label{eq:balmatr}
\fM_{\Xbar} = \left (
  \begin{array}{cccccc}
    -\phi_1(x_1) & \cdots & -\phi_1(x_l) &    \phi_1(y_1) & \cdots & \phi_1(y_m) \\
\vdots & &\vdots &\vdots & & \vdots\\
    - \phi_r(x_1) & \cdots & - \phi_r(x_l)  &    \phi_r(y_1) & \cdots & \phi_r(y_m)\\
  \end{array}
\right )
\end{equation}
Notice that only the orbifold points such that $p_j\neq q_j-1$ comme
into play in this matrix.

Then we define the integers
$$
\fC_1(\Xbar) = \rank \fM_{\Xbar}\quad \mbox{and}\quad \fC_2(\Xbar)
=\dim (C^{k+r}_+\cap\ker\fM_{\Xbar}), 
$$
where $C^{l+m}_+$ is the cone of vectors with positive entries in $\RR^{l+m}$.

We extract from Arezzo-Pacard gluing theory \cite{ArePacX05} the
following theorem. 
\begin{theo}
  \label{theoglue}
Let $\Xbar$ be a compact complex orbifold
surface with isolated singularities $\{x_j,\; 1\leq j\leq k\}$ modelled on
$\CC^2/\Gamma_{p_j,q_j}$  where $(p_j,q_j)$ are coprime integers with
$0<p_j<q_j$, and   let  $\bar\omega$ be an orbifold  CSCK metric on $\Xbar$. 
Consider the minimal resolution 
$\Xhat\to \Xbar$  
and define $\Xtil \to \Xhat$ by performing further (simple) blow-ups at
some smooth points $\{y_1,\cdots,y_r\}\subset \Xbar\setminus \{x_j\}$.
With the above notation, we are assuming moreover that $\fC_1(\Xbar) = \dim_\CC \cV_0(\Xbar)$ and
$\fC_2(\Xbar)\neq 0$.

Then, given $\epsilon>0$, $n\geq 0$, a compact
domain $D\subset \Xbar\setminus \{x_j,y_j\}$, and an a priori
norm $\|\cdot\|$ on $H^2(\Xtil,\RR)$,
there exists 
a CSCK metric  $\omega$ on $\Xtil$ such
that 
$$|\omega-\pi^*\bar \omega|_{C^n(D)}\leq \epsilon\quad \mbox{and}\quad
\|[\omega]-[\pi^*\bar\omega] \|\leq \epsilon.$$
Here, $\pi$ is the canonical map obtained by composition $\Xtil\to \Xhat \to \Xbar$, and $C^n(D)$ is
the norm with $n$ derivatives on the domain $D$ measured w.r.t. the metric $\pi^*\bar \omega$.
\end{theo}
\begin{rmk}
\label{rmk:equi}
  There is also an equivariant version of the above theorem: for
  instance suppose 
  that $\Xbar$ is acted on by a finite group $G$ of holomorphic
  isometries and that the set of points $\{x_j, 1\leq j\leq
  k\}\cup\{y_j, 1\leq j\leq m\}$ is $G$-invariant. If
  there are no non-trivial $G$-invariant holomorphic vector fields,
  then the conclusion of Theorem~\ref{maintheo} holds.
\end{rmk}
Strictly
  speaking,   Arezzo-Pacard  do not   address the orbifold
  case in \cite{ArePacX05}. The proof of their main gluing result  \cite[Proposition
  1.1]{ArePacX05} uses in a crucial way the existence of a $\log$-term
  in the expansion of the Burns metric ( cf. \eqref{eq:asympt}) with
coefficitent  $m>0$. When there are 
  holomorphic vector fields, the gluing problem is obstructed an one
  needs to work orthogonally to the kernel to overcome the difficulty.
The $\log$-term at each $y_j$ is used to do so, and is involved in a
certain balancing 
condition reformulated using the
  matrix~$\fM_{\Xbar}$ (when $\{x_j\}=\emptyset$), and ultimately,
  $\fC_1(\Xbar)$ and $\fC_2(\Xbar)$.

The result extends in a straightforward way to the orbifold case using
the  $\log$-term  in the
expansion of the scalar-flat ALE metric on $Y_{p,q}$ whose coefficient needs to be
computed
 (cf. \S\ref{sec:calsin}), and Theorem \ref{theoglue}
follows. The form of the matrix $\fM_{\Xbar}$
comes from the fact that $m$ has opposite signs for the Burns metric
and the Calderbank-Singer metrics on $Y_{p,q}$ whenever  $p\neq
q-1$. The coefficient is zero in the case $p=q-1$, which is the reason why the
points $x_{l+1},\cdots,x_k$ does not come into play in $\fM_{\Xbar}$.
\begin{rmk}
\label{rmk:laborious}
  It would be much more tedious to deal with the case where
  $\cV_0(\Xbar)\neq 0$ and the matrix $\fM$ is empty, that is when
  $\{y_j\}=\emptyset$ and every orbifold singularity $x_j$
verifies
  $p_j=q_j-1$. Since the ALE metric has no $\log$-term in its expansion,
  one would have to refine all the analysis of \cite{ArePacX04}  and
  work out what the new balancing condition is. It is likely that it
  would require introducting  a new matrix $\fM$ with entries taking into account higher
  derivatives of $\phi_j$.
\end{rmk}

\begin{rmk}
\label{rmk:extension}
A result similar to Theorem \ref{theoglue} holds if one replaces
``CSCK'' with ``SFK''.  The gluing theory for SFK metrics
  of~\cite{RolSin05} holds only when there a no non-trivial
  holomorphic vector fields (similarly to \cite{ArePacX04} for CSCK). However the gluing result can be extended
  in presence of holomorphic vector fields along the same line as 
  \cite{ArePacX05}.
\end{rmk}

\begin{rmk}
There is a very striking interpretation of the matrix
$\fM_{\Xbar}$  in terms of moment map:
the condition $\fC_2(\Xbar)\neq 0$ implies that 
\begin{equation}
  \label{eq:kn}
{\mu}(x_1,\cdots,x_l,y_1,\cdots,y_m) := \fM_{\Xbar}\left (
    \begin{array}{c}
      a_1 \\ \vdots \\a_{m+l}
    \end{array}
\right ) =0
\end{equation}
for some $a_j>0$. Then we form the product $W=\Xbar^{m+l}$ endowed with
the symplectic form $\Omega = \sum_{j=1}^{m+l}a_j\pi^*_j\bar\omega$,
where $\pi_j$ is the $j$-th canonical projection
$\Xbar^{m_l}\to\Xbar$. 
For simplicity, we are assuming that the vector field
$\delb^\sharp\phi_j$ are induced by a torus action $\TT^r$. In
particular, the action is commutative. Hence, one can define a
left-action of $t\in \TT^{r}$ on $(W,\Omega)$ as follows:
$$
t\cdot (z_1,\cdots,z_l,z_{l+1},\cdots,z_{m+l}) =  (t^{-1}\cdot
z_1,\cdots,t^{-1}\cdot z_l,t\cdot z_{l+1},\cdots,t\cdot z_{m+1}).
$$
 Then, $\mathrm{\mu}$ is  a moment map for this action. Thanks to the
Kempf-Ness theorem,  
 the condition \eqref{eq:kn} implies that the configuration of points
 $(x_1,\cdots, 
 x_l,y_1,\cdots, y_l)\in W$ is
polystable in the  GIT sense  (the reader can consult nice surveys
\cite{BiqBou04,Tho06}  for more backgroung material). Moreover, when
$\Xbar$ is smooth, this stability
condition which 
garantees the existence of 
a CSCK metric on the blow-up  is also (almost) necessary by
a result of Stoppa \cite{Sto07}.
Sadly, there is at the moment no clear understanding of this fact
when $\{x_j\}\neq\emptyset$, i.e. when orbifold singularities do occur
in Theorem~\ref{theoglue}, and in particular when negative
signs of the matrix
$\fM_{\Xbar}$ come
into play.
\end{rmk}

\subsubsection{Proof of the main results}
\label{sec:proof}
Consider a representation $\rho$ which satisfies 
 $$ \rho(l_j) =
\pm \left ( 
\begin{array}{cc}
  e^{i\pi/q_j} & 0 \\
0 &   e^{-i\pi/q_j}
\end{array}
\right )\quad \mbox{for all $j$}.
$$
More generally, any representation satisfying this property, up to (a
global) conjugation, is called \emph{sporadic}.
\begin{rmk}
  We have chosen this terminology since the representation $\rho$
  associated via Mehta-Seshadri theorem~\ref{prop:MS}
  to a sporadic parabolically
  polystable ruled surface is automatically sporadic in the above sense.
\end{rmk}

We deduce the  following corollary from Theorem~\ref{theoglue}.
\begin{cor}
\label{corglue}
Let $\Sbar$ be a good orbifold Riemann surface and
$\rho\in \cR(\ovS )$. If we have $\cV_0(\Sbar)=0$ and $\dim_\CC
\cV_0(\Mbar_\rho) =1$, we assume moreover that $\rho$ is non-sporadic.

Then, the minimal resolution $\Mhat_\rho$ of $\Mbar_\rho$ carries a CSCK metric~$\omega$. If $\chiorb(\Sbar)<0$, we may assume moreover that
$\omega$ is scalar-flat.
\end{cor}

\begin{proof}
According to \S\ref{sec:korb}, $\Mbar_\rho$ carries a CSCK
metric. 
If $\Sbar$ has no orbifold points, there is nothing to prove because
$\Mbar_\rho$ is already smooth.

Assume that $\Sbar$ has at least one orbifold point. Since $\Sbar$ is
good, either 
\begin{enumerate}
\item \label{enum:cyclic} it is a quotient of the form $\CP^1/\ZZ_q$ or
\item \label{enum:novf} it has no non-trivial holomorphic vector field.
\end{enumerate}
\medskip

\paragraph{{\bf Case}~\ref{enum:cyclic}}
Assume that $\Sbar=\CP^1/\ZZ_q$. We work $\ZZ_2$-invariantly as in the
proof of Proposition~\ref{prop:vf2}.
There
are no non-trivial $\ZZ_2$-invariant holomorphic vector fields on
$\Mbar_\rho$ and one can work using an equivariant version of gluing
theorem (cf. Remark~\ref{rmk:equi}). 
\medskip

\paragraph{{\bf Case}~\ref{enum:novf}} Assume that $\Sbar$ is a good
orbifold Riemann surface with no non-trivial holomorphic vector field.
Then we are either in  case \ref{enum:0} or
\ref{enum:1} of  Proposition \ref{prop:vf}. 
If $\rho$ does
not fix any 
point in $\CP^1$, the space $\cV_0(\Mbar_\rho)$ is
trivial.  This is the case addressed in \cite{RolSin05} and
\cite{RolSinX04}. The gluing theory of \cite{ArePacX04} applies and we
get a CSCK 
metric on  
$\Mhat_\rho$. If $\chiorb(\Sbar)<0$, we apply \cite[Theorem
4.1.1]{RolSin05} and obtain a scalar-flat K\"ahler metric.

Suppose now that  $\cV_0(\Mbar_\rho)$ is $1$ dimensional, or
equivalently, that $\rho$ fixes exactly two points in $\CP^1$.
 Up to an isometry, we may assume that $\rho$ fixes exactly
the points $[0:1]$ and $[1:0]$. 

We introduce the holomorphic vector field on $\CP^1$ 
given by $\delb^\sharp\phi$, where $\phi$ is given at \eqref{eq:pot}.
Notice that $\phi$ is invariant by rotations about 
the axis going through $[0:1]$ and $[1:0]$.
Hence
 $\phi$ descends to a function on $\Mbar_\rho$ and the kernel of
$\LL_{\Mbar}$ is spanned by $1$ and $\phi$. 
The corresponding matrix  is
$$ 
\fM_{\Mbar_\rho} = (-\phi(x_1),\cdots,-\phi(x_l)),
$$
where $x_1,\cdots,x_l,x_{l+1},\cdots,x_k$  are given by orbifold singularities in $\Mbar_\rho$
over each marked point $P_j\in \Sbar$. By definition of $\phi$, we
have $\phi(x_j)=\pm 1$. The assumption that $\rho$ is not sporadic
implies that $\{x_1,\cdots,x_l\}$ is not empty and that $\phi$ is not
constant on this set.
It follows that 
 $ \fM_{\Mbar_\rho}$ is surjective and there is a vector with positive
 entries in its kernel. Hence Theorem \ref{theoglue} applies and the
corollary is proved for the CSCK case.

Under the assumption that $\chiorb(\Sbar)<0$, we know that
$\Mbar_\rho$ can be endowed with am orbifold SFK metric.  Then rely on
Remark~\ref{rmk:extension} to construct a SFK metric on $\Mhat_\rho$, or we use a
simple trick: we can construct a continuous family of CSCK metrics
$\bar \omega_t$ on $\Mhat_\rho$ (with suitable choices for the
curvatures $\kappa_1$
and $\kappa_2$) such that the scalar curvature of
$\bar\omega_0$ is negative whereas the scalar curvature of
$\bar\omega_1$ is positive. Applying Theorem~\ref{theoglue} as before to the
family, we obtain a continuous family of CSCK metrics $\omega_t$ on
$\Mhat_\rho$. If $\epsilon$ is chosen small enough, the scalar
curvature of $\omega_t$ must change sign, hence $\omega_t$ is
scalar-flat for some~$t$.
\end{proof}

The proof of Theorem~\ref{maintheo} is now immediate.
\begin{proof}[Proof of Theorem~\ref{maintheo}]
This is a direct consequence of Proposition~\ref{prop:MS} and 
Corollary~\ref{corglue}.
\end{proof}

The proof of Proposition~\ref{prop:further} goes along the same lines.
\begin{proof}[Proof of Proposition~\ref{prop:further}]
The proof is completely similar to Theorem~\ref{maintheo}. We just
have check that we can allow further blow-ups at some smooth points
$\{y_j\}$ in Corollary~\ref{corglue} to get a CSCK metric on
$\Mtil_\rho$. 

By assumption, we just have to deal with the Case~\ref{enum:1} in the
proof of Corollary~\ref{corglue}, when $\rho$ fixes exactly two points
of $\CP^1$.

The corresponding matrix now
  has now the form
$$
\fM_{\Mbar_\rho} = (-\phi(x_1),\cdots,-\phi(x_l),\phi(y_1),\cdots,\phi(y_m)).
$$
This matrix contains entries of the form $\pm 1$ with both signs.
Therefore it is surjective and has a vector with
positive entries in its kernel. Hence Theorem~\ref{theoglue} still applies.
\end{proof}

\begin{rmk}
  \label{sec:special}
If we allow trivial parabolic structures in
Proposition~\ref{prop:further}, $\Mbar_\rho$ is smooth. If $\rho$ is
trivial, we have $\Mbar_\rho \simeq \Sigma\times\CP^1$. According to
\cite[\S8, Example 5]{ArePacX04} it is possible to blow up the CSCK metric for special
configurations of points $\{y_j\}$.
If $\rho$ is not trivial, it has to fix exactly two points of $\CP^1$ 
since we are
in the strictly parabolic case.
In this case $\cV_0(\Mbar_\rho)= \ip{\delb^\sharp\phi}$ and the corresponding
matrix is $\fM=(\phi(y_1),\cdots,\phi(y_m))$. As soon as this matrix
contains a positive and a negative entry, we have $\fC_1=\fC_2=1$ and
the gluing theorem applies, hence we can blow-up generic
configurations provided $m\geq 2$. 

Suppose that we allow  the case two-punctured sphere in the assumptions
of Proposition~\ref{prop:further}.  Recall that we are using the
$\ZZ_2$-equivariant version of the gluing theorem, in order to prove
Corollary~\ref{corglue}. Hence the conclusion of
Proposition~\ref{prop:further} holds provided that $\{y_j\}$ is a
$\ZZ_2$-invariant set.
\end{rmk}

We can also complete the study of the particular example highlighted
throughout the paper:
\begin{proof}[Proof of Corollary~\ref{propexample}]
  Using the observation of \S\ref{sec:simpex2} and
  \ref{sec:back} we can apply Theorem~\ref{maintheo} and get a CSCK
  metric $\omega$ on $\Mhat$. The qualitative properties of $\omega$
  are obtain by tracking back the construction of the metric via gluing: every 
  K\"ahler class $\Omega$ on $\Mbar$ can be represented by a CSCK
  metric $\bar \omega$. Then the rest of the proposition follows from
  Theorem~\ref{theoglue}    (or more accurately
  Remark \ref{rmk:equi}).
\end{proof}
\begin{rmk}
  \label{rmk:toric}
Consider the parabolic structure on $\CP^1\times\CP^1$ given in
Remark~\ref{rmk:toric0} pictured in the following diagram,
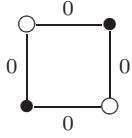
\begin{equation*}
\xymatrix{
 *+[o][F-]{}\ar@{-}[r]^{0} & *{\bullet}{}\\
*{\bullet}{}
\ar@{-}[u]^{ 0} 
\ar@{-}[r]_{0} & *+[o][F-]{}
\ar@{-}[u]_{0} &
}
\end{equation*}
where the black points represent the parabolic points to be blown-up. 
With the notations of \S\ref{secitbup}, the
iterated  blow-up $\Mhat$  for the weights $\alpha_1=\alpha_2$ contains a configuration of
rational curves shown in the diagram below
\begin{equation}
\xymatrix{
*+[o][F-]{}\ar@{-}[r]^{-e_1} &  *+[o][F-]{}
\ar@{--}[r]  &  *+[o][F-]{}
\ar@{-}[r]^{-e_l} &  *+[o][F-]{}
\ar@{-}[r]^{-1} &  *+[o][F-]{}
\ar@{-}[r]^{-e'_{m}} &  *+[o][F-]{}
\ar@{--}[r] &  *+[o][F-]{}
\ar@{-}[r]^{ -e'_{1}} & *+[o][F-]{}\\
*+[o][F-]{}
\ar@{-}[u]^{ -1} 
\ar@{-}[r]_{-e'_1} &  *+[o][F-]{}
\ar@{--}[r]  &  *+[o][F-]{}
\ar@{-}[r]_{-e'_m} &  *+[o][F-]{}
\ar@{-}[r]_{-1} &  *+[o][F-]{}
\ar@{-}[r]_{-e_{l}}  &  *+[o][F-]{}
\ar@{--}[r]  &  *+[o][F-]{}
\ar@{-}[r]_{ -e_{1}} & *+[o][F-]{}
\ar@{-}[u]_{ -1} &
}
\end{equation}
where $e_j$ are the coefficient of the continued fraction expansion of
$\alpha_1=\alpha_2$ and $e'_j$ those for $1-\alpha_1$. 
Similarly to the basic example, $\Mhat$
 admits a
CSCK metric and two linearly independent holomorphic vector
fields. Hence the metric has toric symmetry. 
\end{rmk}

\begin{proof}[Proof of Corollary~\ref{prop:simplest}]
We start with $\Sigma=\CP^1$, $M=\CP^1\times\CP^1$ and the map
$M\to\Sigma$ given by the first projection.
We pick $4$ points, say
$$
P_1=[0:1],\quad P_2=[1:0],\quad P_3=[1:1],\quad P_4=[-1:1]
$$
in $\Sigma$ with weights $\alpha_1=\alpha_2=\frac 12$ and
$\alpha_3=\alpha_4=\frac 13$. We consider two disctinct constant
sections $S_1$ and $S_2$ of $M\to\CP^1$ and declare that $Q_1,Q_3\in S_1$ and
$Q_2,Q_4\in S_2$. It is easy to check that such parabolic structure
make $M$ into a strictly parabolically polystable ruled  surface. It
is moreover non-sporadic since $Q_3\in S_1$ and $Q_4\in S_2$ have both
weight $1/3$.
\end{proof}

\section{LeBrun-Singer's results versus Theorem~\ref{maintheo}}
\label{sec:lsvs}
Strictly parabolically polystable ruled surfaces were studied by
Le\-Brun-Singer more that ten years ago, in relation with  SFK metrics
on their simple blow-ups. As  the proof of
Theorem~\ref{maintheo} was being completed, it became clear that such
result would provide scores of
counterexamples to the classification result deduced from
\cite[Proposition 3.1]{LebSin93} as
explained in \S\ref{sec:counter}. The error
is in fact located   at the end of the proof of that proposition.
Claude LeBrun then pointed out  that the mistake had  already been
spotted in a paper of Kim-LeBrun-Pontecorvo~\cite{KimLebPon97}, 
where the first counter-examples are constructed. 

\subsection{Counterexamples}
\label{sec:counter}
 There are lots of 
 strictly
parabolically polystable ruled surfaces over Riemann surfaces. We 
exhibit two infinite families over the sphere and the torus an let the
reader play this amusing game more generally. 

\subsubsection{Example over the sphere}
Consider The sphere with
$3$-punctures and the 
ruled surface $\pi:\Mch\to\CP^1$ where $\Mch=\CP^1\times \CP^1$ and  $\pi$ is the first
projection. Define the  parabolic structure
with marked points $Q_1 = ([1:0],[1:0])$, $Q_2 = ([0:1],[1:0])$, $Q_3
= ([1:1],[0:1])$ and the weights $\alpha_1=\alpha_2=2/q$ and
$\alpha_3=4/q$ for some odd integer $q\geq 2$.
It is easy to see that $\Mch$ is parabolically polystable but not
stable (proof similar to Lemma~\ref{lemma:ex}). For $q$ large enough,
we have $\chiorb(\Sbar)<0$ and Theorem~\ref{maintheo} provides a
scalar-flat K\"ahler metric on the rational surface $\Mhat\to\CP^1$.
Since $\Mch$ is strictly polystable, we deduce that $\cV_0(\Mhat)$ is
non-trivial, which can be easily seen by hand anyway. 

\subsubsection{Example over the torus}
Consideer the two-punctured  torus $\TT$ with marked points  $P_1$,
$P_2$ and a ruled surface  $\pi:\TT\times\CP^1\to\TT$. Pick two point
$0$ and $\infty$ on $\CP^1$. Let $S_1$ and $S_2$ be the two
constant sections corresponding to $0$ and $\infty$.
 We define $Q_j$ to be
the point over $P_j$ lying on $S_j$. We choose the weights $\alpha_1 =
\alpha_2 = \frac 2q$, where $q$ is odd. By definition we have
$\mu(S_1)=\mu(S_2)=0$ and it is easy to check that $M$ is strictly parabolically
polystable. Morever, we see that the weights are non-sporadic.

Now these examples contradict the
conclusion of  \cite[Proposition 3.1]{LebSin93}  that
the genus of $\Sigma$  should at least $2$.
The the proof can be fixed   (as well as the result deduced
from the proposition in \cite[Theorem 3.7,
Corollary 3.9]{LebSin93})  by adding the extra
condition that the 
$S^1$-action induced by the holomorphic vector field is semi-free.
The reader may refer to 
 \cite[Remark
p. 86]{KimLebPon97}) where this issue  is discussed.
 Notice that the $S^1$-action precisely fails to
be semi-free in our orbifold construction.

\subsection{Comparison with LeBrun's metrics}
\label{sec:compare}
Given a parabolic ruled surface
$\Mch\to\Shat$, we introduce the simple blow-up
at each parabolic point $\Mhat^s\to\Mch$. If $g(\Shat)\geq
2$ and $\Mhat$ is strictly parabolically polystable, LeBrun's
ansatz~\cite{Leb91} (see also \cite[\S3.2]{LebSin93})
gives  a scalar-flat K\"ahler metric $\omega^s$ on $\Mhat^s$ such that  
$$
\frac{[\omega^s]\cdot [E_j]}{[\omega^s]\cdot [F]} = \alpha_j,
$$
where $F$ is a generic fiber, $E_j$ is one of the exceptional curves
and $\alpha_j$ is the
corresponding parabolic weight.

The iterated blow-up $\Mhat$ is obtained by performing further blow-up
on  $\Mhat^s$. Pictorially, we are passing from diagram~\eqref{e2.844}
to diagram~\eqref{eq:mbup}.
 Although the two surfaces are different manifolds, one
can use the blow-up map
$
\pi:\Mhat\to\Mhat^s$ to compare  the K\"ahler classes of each
construction.

Let $E_j'$ be the $-1$ curve corresponding to the last iterated blow-up
of $Q_j$ (cf. \eqref{eq:mbup}). Then $\pi^*[\omega^s]\cdot [E'_j]=0$.
However, we can see from our construction that $[\omega]\cdot [E'_j]
\to [\bar\omega ]\cdot[F_j]\neq 0$ as $\epsilon \to 0$ in the gluing
theorem, where $\bar\omega$ is the scalar-flat K\"ahler metric on
$\Mbar$ and $F_j\subset \Mbar$ is the orbifold fiber over $P_j$.

Therefore, the K\"ahler class $[\omega]$ is in no way a small
perturbation of $\pi^*[\omega^s]$. In particular, expecting to
construct $\omega$ starting from $\omega^s$ and applying further
blow-up and gluing theorem seems hopeless. An interesting open
question is whether 
there is a continuous family of cohomology classes from  $[\omega]$
to $\pi^*[\omega^s]$ representing scalar flat K\"ahler classes on some
blow-up of $\Mch$.

\subsection{Further remarks}
In the minimal case, the correspondence
between CSCK metrics on geometrically ruled surfaces and polystability was first pointed out by Burns-de
Bartolomeis \cite{BurBar88} in the scalar-flat case. The picture was completed by the subsequent
work of LeBrun, in the case of negative scalar curvature~\cite{Leb95},
and, more recently, of
Apostolov-T\o{}nnesen~\cite{ApoTonX04} in full generality, relying on
the ground-breaking work of Donaldson~\cite{Don01}. 

In contrast, 
we seem at the moment quite far away from a complete understanding of how the
stability of bundles is related to CSCK metrics for blown-up ruled
surfaces. More specifically, we do not understand  how stability
of bundles is related to the different notions of geometric stability of
the ruled surface like K-stability.
One of the most exciting result would be to prove a sort of converse
to Theorem~\ref{maintheo}, that is to establish a correspondence
between the two categories:
\vspace{10pt}

\begin{center}
\begin{tabular}{lcr}
 \fbox{\parbox{1.6in}{Parabolically polystable ruled surfaces $\Mch$} }
& $\leftrightarrow$ &\fbox{\parbox{2in}{CSCK metrics on the corresponding
    blow-up $\Mhat$ in ``special'' K\"ahler classes}}
\end{tabular}
\end{center}

\section{Asymptotics of the Joyce-Calderbank-Singer metrics}
\label{sec:asymptotics}
\subsection{Background and notation}

Let $\HH\subset \RR^2$ stand for the half-space with coordinates
$x>0$ and $y$. Let $U$ be an open set of $\HH$ and suppose given a
pair of functions $v_1, v_2:U \to \RR^2$ satisfying
\begin{equation}\label{e1.7.1.7}
\frac{\del v_1}{\del y} = \frac{\del v_2}{\del x},\;\;
x\frac{\del v_1}{\del x} + x\frac{\del v_2}{\del y} = v_1.
\end{equation}
Denote by $\ip{v_1,v_2}$ the determinant of the $2\times 2$ matrix
whose rows are $v_1$ and $v_2$.  The following result appears in the
literature [Joyce, Calderbank--Singer, LeBrun]:
\begin{theo} 
\label{theo:jcs}
Let $U_0\subset U$ be the open set (assumed non-empty) on
  which $\ip{v_1,v_2}\not=0$. Let $M = U_0\times T^2$ with flat
  (multi-valued) coordinates $(t_1,t_2)$ on $T^2$. Then
\begin{equation}\label{e3.14.11.03}
g
=\frac{ x\langle v_1, v_2\rangle}{2(x^2+y^2)}\left(
\frac{\rd x^2 + \rd y^2}{x^2}
+ \frac{\langle v_1, \rd t \rangle^2 +  \langle v_2, \rd t \rangle^2}{
\langle v_1, v_2\rangle^2}\right)
\end{equation}
is a half-conformally flat metric on $M$.  

Let
\begin{equation}\label{e1.13.1.7}
x = r^{-2}\sin 2\theta,\;y = r^{-2}\cos 2\theta
\end{equation}
where $r>0$ and $0\leq \theta \leq \pi/2$, and write $s=\sin\theta, c
= \cos\theta$.  Then the almost-complex
structure $J$ given by
\begin{equation}\label{e2.13.1.7}
 J\rd r = rsc\frac{\ip{v_2,\rd t}}{\ip{v_1,v_2}},\;
 J\rd \theta = -sc\frac{\ip{v_1,\rd t}}{\ip{v_1,v_2}}
\end{equation}
is integrable and $g$ is K\"ahler with respect to $J$, with K\"ahler
form
\begin{equation}\label{e3.13.1.7}
\omega = \left(r\rd r\wedge\ip{v_2,\rd t}
- r^2\rd \theta\wedge\ip{v_1,\rd t}\right).
\end{equation}
In particular $g$ is a scalar-flat K\"ahler metric on $M$.
\label{t1.15.1.7}
\end{theo}
\begin{proof}
We refer to the literature for the proof that $g$ is half-conformally
flat. To show that $J$ is integrable, note that $J^2=-1$ and
\eqref{e2.13.1.7}
imply together that
\begin{equation}\label{eq:cxstruct}
J\rd t = \frac{1}{rsc}\left( v_1\otimes \rd r + v_2\otimes r\rd
  \theta\right).
\end{equation}
(Recall that $t=(t_1,t_2)$, so $J\rd t$ is really a pair of $1$-forms
on $M$.) In particular each component of $J\rd t - i\rd t$ is a
$(1,0)$-form. We claim that $\rd (J\rd t - i \rd t) = 0$, and this
certainly implies that $J$ is integrable. Indeed,
\begin{equation}\label{e4.13.1.7}
\rd (J\rd t - i\rd t) = 
\rd J\rd t = \frac{\rd r\wedge\rd \theta}{rsc}\left(
r\frac{\del v_2}{\del r} - \frac{\del v_1}{\del \theta} - 
\frac{s^2-c^2}{sc}v_1\right).
\end{equation}
But by the chain rule,
\begin{equation}\label{e5.13.1.7}
r\del_r = - 2(x\del_x + y\del_y),\;\;
\del_\theta = - 2(x\del_y - y\del_x);
\end{equation}
substitution of these into \eqref{e4.13.1.7} and use of
\eqref{e1.7.1.7} now proves the claim.

To verify the compatibility of $g$ and $J$, note first that routine
computation gives
\begin{equation}\label{e6.13.1.7}
\frac{\rd x^2 +\rd y^2}{x^2}
=
\frac{\rd r^2 + r^2\rd \theta^2}{r^2s^2c^2}
\end{equation}
while
\begin{equation}\label{e7.13.1.7}
J\rd r^2 + r^2 J\rd \theta^2
= r^2s^2c^2\frac{\ip{v_1,\rd t}^2+\ip{v_2,\rd t}^2}{\ip{v_1,v_2}^2}
\end{equation}
so
\begin{equation}\label{e8.13.1.7}
g = \frac{|\ip{v_1,v_2}|}{sc}\left(
\rd r^2 + J \rd r^2 + r^2(\rd\theta^2 + J\rd\theta^2)\right).
\end{equation}
This shows that $g$ is $J$-hermitian with fundamental $2$-form
\begin{equation}\label{e9.13.1.7}
\omega = \frac{|\ip{v_1,v_2}|}{sc}\left(
\rd r\wedge J\rd r + r^2\rd\theta\wedge J\rd\theta\right).
\end{equation}
Inserting \eqref{e2.13.1.7} here yields \eqref{e3.13.1.7}.

It remains to check that $\rd\omega=0$. We have
\begin{equation}\label{e10.13.1.7}
\rd\omega =
-r \rd r\wedge\rd\theta\wedge
\ip{\del_\theta v_2 - 2v_1 - r\del_r v_1,\rd t}.
\end{equation}
It is now easy to check using \eqref{e5.13.1.7} and \eqref{e1.7.1.7}
that $\rd\omega=0$. The proof is complete.
\end{proof}

\subsection{The flat metric}
Let
\begin{equation}\label{e11.10.1.7}
v_1 = \frac{x}{2}\frac{(1,-1)}{\sqrt{x^2+y^2}},\;\;
v_2 = \frac{1}{2}\left(
  \frac{y(1,-1)}{\sqrt{x^2+y^2}} + (1,1)\right).
\end{equation}
\begin{prop} We have $\ip{v_1,v_2}>0$ in $c\HH$ and $(g,J,\omega)$ is
  the standard flat metric on $\CC^2$.
\end{prop}
\begin{proof} Denote standard linear complex coordinates on $\CC^2$ by
  $(z_1,z_2)$. Introduce coordinates $(r,\theta,\phi,\psi)$ on $\CC^2$
  by setting
\begin{equation}\label{e12.10.1.7}
z_1 = r\cos\theta e^{i\phi},\; z_1 = r\sin\theta e^{i\psi}
\end{equation}
where $\theta \in [0,\pi/2]$ and $\phi,\psi \in [0,2\pi]$. Then the
flat metric $g_0$ becomes
\begin{equation}\label{e13.10.1.7}
g_0 = \rd r^2 + r^2\rd\theta^2 + r^2(c^2\rd \phi^2 + s^2\rd\psi^2)
\end{equation}
and the K\"ahler form is
\begin{align}\nonumber
\omega_0 &= \frac{i}{2}\left(\rd z_1\wedge \rd \overline{z}_1
+ \rd z_2\wedge \rd \overline{z}_2\right)\\
&= 
r\rd r\wedge(c^2\rd\phi + s^2\rd\psi)
- r^2sc\rd\theta\wedge(\rd \phi - \rd\psi). \label{e14.10.1.7}
\end{align}
These coordinates are matched up to the $(x,y,t_1,t_2)$ system by
\begin{equation}\label{e15.10.1.7}
x = r^{-2}\sin 2\theta, y = r^{-2}\cos 2\theta, t_1 = -\psi, t_2 = \phi.
\end{equation}
This is a matter of simple computation. For example
\begin{equation}
v_1 = sc(1,-1),\; v_2 = (c^2, s^2), \ip{v_1,v_2} = sc
\end{equation}
so
\begin{equation}
\ip{v_1,\rd t} = sc(\rd t_1 +\rd t_2),\;
\ip{v_2,\rd t}= -s^2\rd t_1 +c^2\rd t_2.
\end{equation}
Hence it is now easy to verify that in this case the metric \eqref{e3.14.11.03}
agrees with the flat metric $g_0$ if $t_1 = \pm \psi$ and $t_2 = \pm
\phi$. In this case the formula \eqref{e3.13.1.7} for the K\"ahler
form becomes
\begin{equation}
\omega = r\rd r\wedge(-s^2\rd t_1 + c^2\rd t_2) - r^2sc\rd\theta
\wedge(\rd t_1 + \rd t_2).
\end{equation}
Matching this up with the formula for $\omega_0$ fixes the
identification of the angular variables as claimed in
\eqref{e15.10.1.7}.
\end{proof}

We can also see that $r^2/4$ is the K\"ahler potential for this flat
metric, as it should be:
\begin{eqnarray}
\rd J\rd r^2/2 &=& \rd (r/2) J\rd r \\
& = & \rd \left(r^2sc\frac{\ip{v_2,\rd t}}{2\ip{v_1,v_2}}\right)\\
&=& \rd (r^2/2)(c^2\rd\phi + s^2\rd\psi) \\
&=& r\rd r\wedge(c^2\rd\phi + s^2\rd\psi) -
r^2sc\rd\theta\wedge(\rd\phi-\rd\psi) \\
&=& \omega_0.
\end{eqnarray}
\label{s1.15.1.7}
\subsection{More general ALE metrics}

In \cite{CalSin04} more general SFK metrics were considered, where
$(v_1,v_2)$ is a finite linear combination of basic solutions
\begin{equation}
\left(\frac{x}{\sqrt{x^2+(y-a)^2}},\frac{y-a}{\sqrt{x^2+(y-a)^2}}\right).
\end{equation} 
More precisely,
pick a strictly decreasing sequence $\infty\geq y_0 > y_1 >\ldots >
y_{k+1}=0$ and a sequence of pairs $(a_j,b_j)\in \ZZ^2$. Define
\begin{align}
v_1 &= \frac{x}{2}\left(
\frac{1}{\sqrt{x^2+y^2}}(a_{k+1},b_{k+1})
+ \sum_{j=0}^k\frac{1}{\sqrt{x^2 + (y-y_j)^2}}(a_j,b_j)
\right)\label{e1.15.1.7}\\
v_2 &= \frac{1}{2}\left(
\frac{y}{\sqrt{x^2+y^2}}(a_{k+1},b_{k+1})
+ \sum_{j=0}^k\frac{y-y_j}{\sqrt{x^2 + (y-y_j)^2}}(a_j,b_j)
\right)\label{e2.15.1.7}
\end{align}
Then $(v_1,v_2)$ satisfies \eqref{e1.7.1.7}.  If $y_0=\infty$, then
the $j=0$ contribution to $v_1$ is zero, and the contribution to $v_2$
is $-(a_0,b_0)$. We recover the flat metric by taking $k=0$ and
$y_0=\infty$.

From now on, we shall assume that
\begin{equation}
(a_j,b_j) = (m_j-m_{j+1}, n_j- n_{j+1})
\end{equation}
where $(m_j,n_j)$ is a sequence of integer pairs with 
\begin{equation}
(m_0,n_0) = -(m_{k+2},n_{k+2})= (0,-1), (m_1,n_1) = (1,0)
\end{equation}
and
\begin{equation}
m_j>0 \mbox{ for all }j=1,\ldots,k+1;\;\;
m_jn_{j+1} - m_{j+1}n_j = 1\mbox{ for all }j.
\end{equation}
From \cite{CalSin04}, we know that these conditions imply that
$\ip{v_1,v_2}>0$ in $\HH$, so that $g$ is defined on the whole of
$M_0 = \HH\times T^2$. Furthermore, $g$ extends as a smooth metric
to a partial compactification $M$ of $M_0$.  The space $M$ can be
defined in the following way. Let $\overline{M} = \HH\times T^2\cup S^1\times
T^2$ be the manifold with boundary obtained by replacing $\HH$ by
its conformal compactification $\oH = \HH \cup
\{(0,y):y\in\RR\}\cup \infty$. The space $M$ is obtained by blowing
down a circle over the interior of each interval $(y_j,y_{j-1})$;
by blowing down the whole $T^2$ over each of the $y_j$, $j=0,\ldots,
k+1$; and finally by deleting the point corresponding to $y_{k+1}$.
Thus $M$ is non-compact, with its asymptotic region corresponding to a
neighbourhood of $(0,0)$ in $\oH$.

\subsection{Hirzebruch-Jung resolutions}
\label{sec:hj}
The most important application of this construction arises in the
following way.  Let $(p,q)$ be a coprime pair of integers with
$0<p<q$. Then we may consider the orbifold 
$\CC^2/\Gamma_{p,q}$ (recall that the action of $\Gamma_{p,q}$ is
generated by~\eqref{eq:action}).
 The minimal resolution $Y_{p,q}\to\CC^2/\Gamma_{p,q}$  has a toric description by taking
the (modified) continued-fraction expansion~\eqref{e1.844}
and defining $(m_j,n_j)$, for $j\geq 2$, as the $j$-th approximant to
$p/q$:
\begin{equation}
\frac{n_{j+1}}{m_{j+1}}=
\cfrac{1}{e_1-\cfrac{1}{e_2- \cdots \cfrac{1}{e_j}}}.
\end{equation}

\begin{theo} 
\label{theo:asymptmichael}
 Let $(v_1,v_2)$ be defined by 
\eqref{e1.15.1.7} and \eqref{e2.15.1.7} and let the corresponding
structures of Theorem~\ref{t1.15.1.7} be denoted by
$(g_v, J_v, \omega_v)$.  Then this triple defines an asymptotically
locally euclidean structure on  $M$, and 
an approximate K\"ahler potential for $\omega_v$ is given by
\begin{equation}
\label{eq:formf}
\frac{f}{q} = \frac{r^2}{4} + \frac{a+b}{2}\log r + \frac{a-b}{2}c^2 ,
\end{equation}
that is, $ \omega_v = \rd J_v\rd f +\cO(r^{-4})$. Here
$a$ and $b$ are defined by
\begin{equation}\label{e1.8.1.7}
a(q,p) - b(0,1) = \sum_{j=0}^k y_j^{-1}(a_j,b_j).
\end{equation}
\end{theo}
\begin{proof} 
It was explained in \cite{CalSin04} and \cite{Joy95} that $g_v$
  extends to $M$.  To verify the statement about the asymptotics, we
  must expand $(v_1,v_2)$ for small $r$. We have
\begin{equation}\label{e6.14.11.03}
v_1 = s c(1 + a r^{-2})(q,p) -sc(1 + br^{-2})(0,1)
O(r^{-4}),
\end{equation}
and
\begin{equation} \label{e3.24.11.03}
v_2 = c^2(q,p) + s^2(0,1) + O(r^{-4}).
\end{equation}
Now define new angular variables $(\phi,\psi)$ so that
\begin{align}
\ip{v_1,\rd t} =& sc(1+ar^{-2})\rd\phi
-sc(1+br^{-2})\rd\psi \label{eq:phipsi} \\
\ip{v_2,\rd t} =& c^2\rd\phi + s^2 \rd \psi
+O(r^{-4}) \nonumber
\end{align}
so as to have
agreement, to leading order, with the flat case, compare
\S\ref{s1.15.1.7}.  The 
fact that the determinant of the transformation from $(t_1,t_2)$ to
$(\phi,\psi)$ is $q$ means that $(\phi,\psi)$
really live on a $q$-fold cover.  Then
\begin{equation}
\ip{v_1,v_2}=
qsc(1 + (as^2+bc^2)/r^2)
+ O(r^{-4}).
\end{equation}
Having made these substitutions, it is clear that $g_v$ differs from
$g_0$ by terms of order $r^{-2}$ and that there will be similar
estimates on the derivatives of $g_v$, as required for an ALE metric. 
The K\"ahler form is
\begin{align}\label{e3.8.1.7}
\omega_v = &
r\rd r\wedge(c^2\theta\rd\phi + s^2\theta \rd \psi)
-
r^2sc\rd \theta((1+a r^{-2})\rd\phi \nonumber\\
&- (1+b r^{-2})\rd\psi) + O(r^{-4})\nonumber \\
= & \omega_0  - sc\rd\theta\wedge(a\rd\phi - b\rd\psi) + O(r^{-4}).
\end{align}
Our first approximation to the K\"ahler potential is the flat
potential $r^2$, but in applying $\rd J\rd$ we have to be careful to
use $J_v$ rather than $J_0$. We have
\begin{align}
\rd J \rd (qr^2/4) =& 
\rd (qr/2)J\rd r \\
=& \rd (r^2/2)(1 - (as^2+bc^2)r^{-2})(c^2\rd\phi + s^2\rd\psi)) \\
=& \omega_0 - (a-b)sc\rd\theta\wedge(c^2\rd\phi + s^2\rd\psi)\nonumber \\
&-(as^2+bc^2)sc\rd\theta\wedge(\rd\psi - \rd\phi).
\end{align}
Hence
\begin{equation}
\rd J_v \rd (qr^2/4) =
\omega_v + 2sc(as^2+bc^2)\rd\theta\wedge(\rd\phi - \rd\psi) + O(r^{-4})
\end{equation}
Now we easily compute
\begin{equation}
\rd J\rd q\log r =- 2sc\rd\theta\wedge(\rd\phi - \rd\psi)
\end{equation}
and
\begin{equation}
\rd J \rd qc^2 = -2sc(s^2-c^2)\rd\theta\wedge(\rd\phi - \rd\psi).
\end{equation}
Combining the last three equations with the identity
$$
2(as^2 +bc^2) = (a+b) + (a-b)(s^2-c^2)
$$
 completes the proof.
\end{proof}
\medskip

%
%

Since $J_v$ describes the complex structure of the Hirzebruch-Jung
resolution $Y_{p,q}\to\CC^2/\Gamma_{p,q}$ we can write the the
potential $f$ on the chart $(\CC^2\setminus 0)/\Gamma_{p,q}$. Let
$z=(z_1,z_2)$ be the standard holomorphic coordinates on $\CC^2$ and
$|z|=\sqrt{|z_0|^2+|z_1|^2}$.
We deduce the following corollary from
Theorem~\ref{theo:asymptmichael}.
\begin{cor}
  \label{cor:asympthol}
In the holomorphic chart $(\CC^2\setminus 0)/\Gamma_{p,q}$, the metric
is given by $\omega_v = dd^c f$, where 
\begin{equation}
\frac{f}{q} = \frac{|z|^2}{4} + \frac{a+b}{2}\log |z| +\cO(|z|^{-1}),
\end{equation}
and $a$ and $b$ are defined by~\eqref{e1.8.1.7}.
\end{cor}
\begin{proof}
We start by introducing the approximate holomorphic coordinates
$\tilde z =(\tilde z_1,\tilde z_2)$ given
by
$$
\tilde z_1 =r\cos\theta e^{i\phi},\quad \tilde z_2 =r\sin\theta e^{i\psi}.
$$
The holomophic $(1,0)$-form $dt-iJdt$ (cf. proof of
Theorem~\ref{theo:jcs}) is in fact given by a pair of holomorphic
$(1,0)$-forms $\gamma_1$ and $\gamma_2$. Moreover, a direct computation
using~\eqref{eq:cxstruct}, \eqref{e6.14.11.03}, \eqref{e3.24.11.03}
and \eqref{eq:phipsi} shows that 
$$
\gamma_j = \frac{d\tilde z_j}{\tilde z_j} + F_j(\tilde z_1,\tilde z_2),
$$
where $F_j$ decays as  $\cO(r^{-3})$. Now put $f_j = - \int_{\tilde
  z}^\infty F_j$ where the integral is taken along the path $s\mapsto
s\tilde z$ from $1$ to $\infty$. Then $f_j = \cO(r^{-2})$ and we
introduce the holomorphic coordinates $z_j:=\tilde z_j\exp(f_j)$. It
follows that 
\begin{equation}
  \label{eq:holcoord}
z_j = \tilde z_j (1+\cO(r^{-2})).
\end{equation} Notice that the above
computations correct a serie of typos in \cite[Section 5.1]{RolSin05}.

In particular it follows from~\eqref{eq:holcoord} that $|z|\simeq r$
hence by Theorem~\ref{theo:asymptmichael} 
$\omega_v = dd^c f + \beta$ where $\beta$ is a closed $(1,1)$-form
with decay $\cO(|z|^{-4})$. Now the $dd^c$-lemma~\cite[Theorem 8.4.4]{Joy00} shows that $\beta=dd^cu$ for
some function $u=\cO(|z|^{-2})$. Eventually  $\omega_v = dd^c h$ on
$(\CC^2\setminus 0)/\Gamma_{p,q}$ for
$h= f + u$. 

Using~\eqref{eq:holcoord} again with~\eqref{eq:formf}, we see that
$$\frac hq = \frac{|z|^2}4 + \frac{a+b}2\log |z| +\cO(1).
 $$
Using
the fact that the metric is SFK, the
bootstrapping argument given by \cite[Lemma 7.2]{ArePacX04} shows that
we have actually
$$
\frac hq = \frac{|z|^2}4 + c_1 \log |z| + c_2 |z| +c_3 + \cO(|z|^{-1}).
$$
By uniqueness of the expansion, we have $c_1= \frac{a+b}2$.  Since
$dd^c(c_2 |z| +c_3)=0$ we can always assume that the potential of the
K\"ahler form is given by 
$$\frac hq = \frac{|z|^2}4 + \frac{a+b}2\log |z| +\cO(|z|^{-1}),
 $$
and the corollary is proved.
\end{proof}

\subsection{The sign of the log-term}
\label{sec:sign}
It was observed by LeBrun that the coefficient of the log-term is
positive for the Burns metric, zero for the Eguchi-Hanson metric
(which corresponds to $p/q =1/2$) and negative for the spaces
corresponding to $p/q =1/q$ for $q>2$.

More generally, we have
\begin{theo} \label{theo:sign}
Let $M$ be the minimal resolution corresponding to the
  fraction $p/q$. Then the coefficient of the log-term is non-positive
  and is zero if and only if $p = q -1$ (the case that our SFK metric
  is K\"ahler-Einstein).
\end{theo}
\begin{proof}
We have to understand $\mu =a +b$, where the $(m_j,n_j)$ arise from the
continued-fraction expansion of $p/q$ as above and
\begin{equation}\label{e2.14.1.7}
a(q,p) + b(0,-1) = \sum_{j=0}^k y_j^{-1}(m_j-m_{j+1}, n_j-n_{j+1}).
\end{equation}
Write $c_j = y_j^{-1}$, so that $c_0 < c_1 < \ldots < c_k$. From
\eqref{e2.14.1.7}, we find
\begin{eqnarray}
qa &=& \sum_{j=0}^{k+1}(c_j-c_{j-1})m_j\nonumber \\
qb &=& \sum_{j=0}^{k+1}(c_j-c_{j-1})(pm_j-qn_j)
\end{eqnarray}
where by convention $c_{-1} = c_{k+1}=0$. Hence
\begin{align*}
\mu =& c_0 +(c_1-c_0)(m_1 + pm_1 - qn_1)+
(c_2-c_1)(m_2 + pm_2 - qn_2) \\
&+\ldots+ (c_k-c_{k-1})(m_k + pm_k - qn_k) - c_k.
\end{align*}
The $m_j$ are positive, so if we introduce $u_j = m_j(c_j-c_{j-1})$,
then the $u_j$ are positive, 
$$
c_k - c_0 = u_1/m_1 + u_2/m_2 + \cdots + u_k/m_k
$$
so we obtain the formula
\begin{equation}\label{e3.14.1.7}
\mu = 
\sum_{j=1}^k\left(\frac{p}{q} -\frac{n_j}{m_j}
+ \frac{1}{q} -\frac{1}{m_j}\right)u_j
\end{equation}
Thus our result will follow from the
\begin{lemma} Let $p/q<1$ and the $n_j/m_j$ be the continued-fraction
  approximants as before. Then
$$
\frac{p}{q} -\frac{n_j}{m_j} +
\frac{1}{q} -\frac{1}{m_j} \leq 0
$$
for all $j$ with equality if and only if $e_1=e_2=\cdots = e_{k}=2$.
\end{lemma}
\begin{proof}  Let
$$
\delta_j = \frac{n_{j+1}}{m_{j+1}} - \frac{n_j}{m_j} +
\frac{1}{m_{j+1}} - \frac{1}{m_{j}}.
$$
Since
$$
\delta_j = -\frac{m_{j+1}-m_j -1}{m_jm_{j+1}}
$$
and the $m_j$ are strictly increasing, we see that $\delta_j\leq 0$
with equality if and only if $m_{j+1}-m_j=1$.  But 
$$
m_{j+1} = e_jm_j - m_{j-1}
$$
so
$$
m_{j+1}-m_j = (e_j-1)m_j - m_{j-1} \geq m_{j} -m_{j-1}
$$
with equality if and only if $e_j=2$. Hence $m_{j+1} - m_j\geq 
m_1-m_0=1$ with equality if and only if $e_1=\ldots= e_j=2$.  On the
other hand, since $(q,p)=(m_{k+1},n_{k+1})$,
$$
\frac{p}{q} -\frac{n_j}{m_j} +
\frac{1}{q} -\frac{1}{m_j} = \delta_j+\cdots+\delta_{k}\leq 0
$$
with equality if and only if $\delta_j=\ldots=\delta_k=0$ if and only
if $e_1=\ldots = e_k=2$.
\end{proof}
The theorem follows as stated, given the simple observation that the
continued fraction with $e_1=\cdots=e_k=2$ gives $p/q = k/(k+1)$.
\end{proof}

\subsection{Example 1}

The examples considered by LeBrun correspond to the data
$$
(m_0,n_0) = (0,-1),\; (m_1,n_1) = (1,0),\;(m_2,n_2)=(q,1).
$$
In this case, there is only one term in the expression for $\mu$ and we have
$$
\mu = \left(\frac{2}{q} - 1\right)u_1.
$$
Since the Burns metric corresponds to the case $q=1$, we see that
$\mu>0$ in this case, but $\mu=0$ for Eguchi-Hanson ($q=2$) and
$\mu<0$ for all $q>2$.

\subsection{Example 2}

If $y$ is the common endpoint of intervals $I$ labelled by $(m,n)$ and $I'$
labelled by $(m',n')$, then the blow-up at $y$ is obtained by inserting an
additional interval $I''$ between $I$ and $I'$ and giving it the label
$(m+m',n+n')$. Obviously this will destroy the monotonicity the
sequence $n_j/m_j$, but the derivation of \eqref{e3.14.1.7} did not
depend upon this and remains valid. In this case we do not get a sign
for the coefficient of the log-term in general. Indeed the blow-up at
$y_1$ of the previous example corresponds to the data
\begin{align*}
(m_0,n_0) = (0,-1),& \; (m_1,n_1) = (1,0),\\
(m_2,n_2)=(q+1,1),& \; (m_3,n_3)=(q,1).
\end{align*}
From \eqref{e3.14.1.7},
$$
\mu = \left(\frac{2}{q}-1\right)u_1 + 
\left(\frac{2}{q}-\frac{2}{q+1}\right)u_2.
$$
The coefficient of $u_2$ is positive for all $q\geq 1$, so $\mu>0$ if
$q=1$, is positive if $q=2$, and can have either sign, depending upon
the values chosen for $u_1$ and $u_2$, if $q>2$.


\subsection{Mass}
  Let $(Y,g)$ be a Riemannian manifold of 
  dimension $4$ such that there are compact
  set $K\subset Y$ and $K'\subset \CC^2/\Gamma_{p,q}$ and a
  diffeomorphism $\Phi:\CC^2/\Gamma_{p,q}\setminus K' \rightarrow
  Y\setminus K$ such that the metric is given on $Y\setminus
  K$ by 
\begin{equation}
\label{eq:dec}
 g =   \geuc +\cO(|z|^{-\tau}),
\end{equation}
for some $\tau>0$. 

Recall that the mass of $g$ is then given by
$$
m(g) = \frac 1{4\pi^2}\lim_{R\to \infty} \int_{S_R} (\del_i
g_{ij}-\del_jg_{ii})\nu^j \; dvol_{S_r}
$$
where   $S_R\subset \CC^2$ is a Euclidean sphere of radius $|z|=R$ in
the chart at infinity and 
$\nu$ a unit outer normal. The mass is actually independent of the
choice of the diffeomorphism $\Phi$ provided $\tau>0$ (the reader is
refered to \cite{HerBou98} for a leisury survey of these notions).

According to Corollary~\ref{cor:asympthol}, the SFK metric $g$ on
$Y_{p,q}$ constructed
at \S\ref{sec:hj} can be written (after a suitable scaling)
$$
\omega =dd^c(|z|^2 + m\log|z|^2 + u)
$$
where $u=\cO(|z|^{-1})$.
Hence it has a  well defined mass $m(g)$.
 Moreover, we have the following lemma.
\begin{lemma}
\label{lemma:massm}
With the above notations, the SFK metric defined on $Y_{p,q}$ verifies  $m(g)=m$.
\end{lemma}
\begin{proof}
  By considering the order of decay of the terms in the expansion, we
  see that $ m(g) = cm$ for some constant $c$. We can compute
  the constant using an example of scalar-flat K\"ahler ALE space. For
  instance the Burns metric with potential $|z|^2 + m \log
  |z|^2$. LeBrun computes its mass \cite{Leb88} and we have $m (g) = m$. Therefore 
  $c=1$ and the lemma is proved.   
\end{proof}
The proof of Theorem~\ref{theo:mass} follows:
\begin{proof}[Proof of Theorem~\ref{theo:mass}]
  It is an immediate consequence of Lemma~\ref{lemma:massm} and Theorem~\ref{theo:sign}.
\end{proof}

\vspace{10pt}
\bibliographystyle{habbrv}
\bibliography{$HOME/tex/biblio,$HOME/tex/rollin}

\end{document}